\documentclass[12pt]{article}
\usepackage{latexsym,amsfonts,amssymb}
\usepackage[dvips]{color}
\usepackage{colordvi,multicol}

\def \cal{\mathcal}

\newtheorem{thm}{Theorem}[section]
\newtheorem{cor}[thm]{Corollary}
\newtheorem{lem}[thm]{Lemma}
\newtheorem{pro}[thm]{Proposition}
\newtheorem{defi}[thm]{Definition}
\newtheorem{rem}[thm]{Remark}
\newtheorem{exa}[thm]{Example}

\begin{document}
\title{\bf  Hunt's Hypothesis (H) for  the Sum of Two Independent L\'{e}vy Processes}
\author{}

\maketitle
\date{}

 \centerline{Ze-Chun Hu} \centerline{\small College
of Mathematics, Sichuan University, Chengdu, 610064, China}
\centerline{\small E-mail: zchu@scu.edu.cn}
\vskip 0.7cm \centerline{Wei Sun} \centerline{\small Department of
Mathematics and Statistics, Concordia University,}
\centerline{\small Montreal, H3G 1M8, Canada} \centerline{\small
E-mail: wei.sun@concordia.ca}


\vskip 1cm


\vskip 0.5cm \noindent{\bf Abstract}\quad Which L\'{e}vy processes
satisfy Hunt's hypothesis (H) is a long-standing open problem in
probabilistic potential theory. The study of this problem for
one-dimensional L\'{e}vy processes suggests us to consider (H)
from the point of view of the sum of L\'{e}vy processes. In this
paper, we present theorems and examples on the validity of (H) for
the sum of two independent L\'{e}vy processes. We also give a
novel condition on the L\'evy measure which implies (H) for a
large class of one-dimensional L\'{e}vy processes.

\smallskip


\smallskip


\tableofcontents

\section{Introduction}

Let $X$ be a time-homogeneous Markov process. Hunt's hypothesis
(H) says that ``every semipolar set of $X$ is polar". This
hypothesis plays a crucial role in probabilistic potential theory.
In particular, it is equivalent to many important principles of
potential theory under mild conditions. These include the bounded
positivity principle, bounded energy principle, bounded maximum
principle and the bounded regularity principle (see e.g.
\cite[Proposition 1.1]{HSZ15}).

In spite of its importance, (H) has been verified only in special
situations.  About fifty years ago, Professor R.K. Getoor
conjectured that essentially all L\'{e}vy processes satisfy (H).
This conjecture stills remains open and is a major unsolved
problem in the potential theory for L\'{e}vy processes (cf.
\cite[page 70]{B96}).

In the following, we will use a diagram to summarize some
sufficient conditions that obtained so far for the validity of (H)
for L\'{e}vy processes. Let $(\Omega,{\cal F},P)$ be a probability
space and $X=(X_t)_{t\ge 0}$ be an $\mathbf{R}^n$-valued L\'{e}vy
process on $(\Omega,{\cal F},P)$ with L\'{e}vy-Khintchine exponent
$\psi$, i.e.,
\begin{eqnarray*}
E[\exp\{i\langle z,X_t\rangle\}]=\exp\{-t\psi(z)\},\  z\in
\mathbf{R}^n,\ t\ge 0.
\end{eqnarray*}
Hereafter $E$ denotes the expectation w.r.t. (with respect to)
$P$, $\langle\cdot,\cdot\rangle$ and $|\cdot|$ denote respectively
the Euclidean inner product and norm of $\mathbf{R}^n$. The
classical L\'{e}vy-Khintchine formula tells us that
\begin{eqnarray*}
\psi(z)=i\langle a,z\rangle+\frac{1}{2}\langle
z,Qz\rangle+\int_{\mathbf{R}^n} \left(1-e^{i\langle
z,x\rangle}+i\langle z,x\rangle 1_{\{|x|<1\}}\right)\mu(dx),
\end{eqnarray*}
where $a\in \mathbf{R}^n,Q$ is a symmetric nonnegative definite
$n\times n$ matrix, and $\mu$ is a measure (called the L\'evy
measure) on $\mathbf{R}^n\backslash\{0\}$ satisfying
$\int_{\mathbf{R}^n\backslash\{0\}} (1\wedge
|x|^2)\mu(dx)<\infty$.

We use Re$(\psi)$ and Im$(\psi)$ to denote respectively the real
and imaginary parts of $\psi$, and use also $(a,Q,\mu)$ to denote
$\psi$. Define
$$
A:=1+{\rm Re}(\psi),\ \ B:=|1+\psi|.
$$
For a finite (positive) measure $\nu$ on $\mathbf{R}^n$, we denote
$$ \hat\nu(z):=\int_{\mathbf{R}^n}e^{i\langle z,x\rangle}\nu(dx).
$$
$\nu$ is said to have finite 1-energy if
$$
\int_{\mathbf{R}^n}\frac{A(z)}{B^{2}(z)}|\hat \nu(z)|^2dz<\infty.
$$
Throughout this paper, we use $\log$ to denote $\log_e$.

We state below the various sufficient conditions for the validity
of (H) for L\'{e}vy processes.

(ND): $Q$ is non-degenerate, i.e., the rank of $Q$ equals $n$.

(KF): $X$ has resolvent densities w.r.t. the Lebesgue measure and
the Kanda-Forst condition holds, i.e., $|\mbox{Im} (\psi)|\leq cA$
for some  constant $c>0$.

(R): $X$ has resolvent densities w.r.t. the Lebesgue measure and
Rao's condition holds, i.e., $|{\rm Im}(\psi)|\leq Af(A)$, where
$f$ is a positive increasing function on $[1,\infty)$ such that
$\int_N^{\infty}(\lambda f(\lambda))^{-1}d\lambda=\infty$ for some
$N\geq 1$.

(EKFR): $X$ has resolvent densities w.r.t. the Lebesgue measure
and the following extended Kanda-Forst-Rao condition holds:

\noindent There are two measurable functions $\phi_{1}$ and
$\phi_{2}$ on $\mathbf{R}^n$ such that
$\rm{Im}\psi=\phi_{1}+\phi_{2}$, and
$$
|\phi_{1}|\leq Af(A),\ \
\int_{{\mathbf{R}^n}}\frac{|\phi_{2}(z)|}{B^2(z)}dz<\infty,
$$
where  $f$ is a positive increasing function on $[1,\infty)$ such
that $ \int_N^{\infty}(\lambda f(\lambda))^{-1}d\lambda=\infty $
for some $N\ge 1$.

$(C^{B/A})$: $X$ has resolvent densities w.r.t. the Lebesgue
measure and there\ exists a constant $c>0$ such that $B(z)\le
cA(z)\log (2+B(z))[\log\log(2+B(z))],\ \forall z\in \mathbf{R}^n.$

$(C^0)$: $X$ has resolvent densities w.r.t. the Lebesgue measure
and for any finite measure $\nu$  on $\mathbf{R}^n$  of finite
1-energy,
$$
\int_{\mathbf{R}^n}
\frac{1}{B(z)\log(2+B(z))[\log\log(2+B(z))]}|\hat
\nu(z)|^2dz<\infty.
$$

(SYM): $X$ has resolvent densities w.r.t. the Lebesgue measure and
is symmetric.

(SP): $X$ has bounded continuous transition densities,
 and $X$ and its symmetrization  have the same polar sets.

(S): $\mu({\mathbf{R}^n\backslash \sqrt{Q}\mathbf{R}^n})<\infty$
and  the following solution condition holds:

The equation $ \sqrt{Q}y=-a-\int_{\mathbf{R}^n\backslash
\sqrt{Q}\mathbf{R}^n}x1_{\{|x|<1\}}\mu(dx) $ has at least one
solution $y\in \mathbf{R}^n$.

Now we can present  the diagram that summarizes all the above sufficient
conditions for the validity of (H) for L\'{e}vy processes.
\begin{eqnarray*}
\begin{array}{ccccccccc}
& & (SYM)&&&&&&\\
& & \Downarrow &&&&&&\\
(ND)&\Rightarrow& (KF)& \Rightarrow & (C^{B/A})&\Rightarrow &(C^0)&&\\
&&\Downarrow&&&&\Downarrow&&\\
&&(R)& \Rightarrow & (EKFR)& \Rightarrow& (H)&\Leftarrow& (SP) \\
& & &&&&\Uparrow&&\\
& & &&&&(S)&&
\end{array}
\end{eqnarray*}
We refer the readers to \cite{Ka76,F75,R88,HS12,HSZ15,HS16} for
the proof of the diagram. We also refer the readers to \cite{HN16}
and \cite{Fi14} for recent interesting results on the validity of
(H). In \cite{HN16}, Hansen and Netuka showed that (H) holds if
there exists a Green function $G>0$ which locally satisfies the
triangle inequality $G(x,z)\wedge G(y,z)\leq CG(x,y)$. In
\cite{Fi14}, Fitzsimmons  showed that Gross's Bwownian motion,
which is an infinite-dimensional L\'evy process, fails to satisfy
(H).

In this paper, we will further study Hunt's hypothesis (H) from
the point of view of the sum of two independent L\'{e}vy
processes. The rest of the paper is organized as follows. In
Section 2, we discuss (H) for one-dimensional L\'{e}vy processes
and provide motivation for exploring (H) through considering sums
of L\'{e}vy processes. Theorem \ref{thm2.5} below extends a result
of Kesten \cite{Ke69}, and Theorem \ref{thm-2.6} below presents a
novel condition on the L\'evy measure $\mu$ which implies (H) for
a large class of one-dimensional L\'{e}vy processes. In Section 3,
we consider (H) for the sum of two independent L\'{e}vy processes
without assuming that resolvent densities exist. We show that if
$X_1$ satisfies (H) and $X_2$ is a compound Poisson process, then
$X_1+X_2$ satisfies (H); and that if both $X_1$ and $X_2$ satisfy
condition (S), then $X_1+X_2$ satisfies (H). In Section 4, we
consider (H) for the sum of two independent L\'{e}vy processes
under the assumption that resolvent densities exist. Roughly
speaking, the results imply that if $X_1$ satisfies (H) and $X_2$
is suitably controlled by $X_1$, then $X_1+X_2$ satisfies (H).

\section{(H) for one-dimensional L\'{e}vy processes}\setcounter{equation}{0}
In this section, we consider Hunt's hypothesis (H) for
one-dimensional L\'{e}vy processes. Let $X=(X_t)_{t\geq 0}$ be a
L\'{e}vy process on $\mathbf{R}$ with L\'{e}vy-Khintchine exponent
$\psi$ and $(a,Q,\mu)$, where $Q$ is a nonnegative constant. If
$\int (1\wedge|x| )\mu(dx)<\infty$, we write
\begin{eqnarray*}
\psi(z)=ia'z+\frac{1}{2}Q z^2+\int_{\mathbf{R}}
\left(1-e^{i\langle z,x\rangle}\right)\mu(dx).
\end{eqnarray*}

\subsection{Motivation}

Let us start by recalling a beautiful result of Bretagnolle
\cite{Br71}. Define
\begin{eqnarray}\label{sec2-a}
{\cal C}=\{x\in \mathbf{R}: P\{X_t=x\ \mbox{for some}\ t>0\}>0\},
\end{eqnarray}
and consider the following different cases:
\begin{itemize}
\item[A.] $Q>0$.

\item[B.] $Q=0; \int (1\wedge |x|)\mu(dx)=+\infty$.

\item[C.] $Q=0; \int (1\wedge |x|)\mu(dx)<+\infty$. We further
decompose it into the following three subcases:
\begin{itemize}
\item[$C_1.$] $a'=0$, \item[$C_2.$] $a'>0$, $\mu$ does not charge
$\mathbf{R}^-:=\{x\in \mathbf{R}:x<0\}$. \item[$C_3.$] $a'>0$,
$\mu$ charges $\mathbf{R}^-$.
\end{itemize}
\end{itemize}

\begin{thm} \label{thm1}(Bretagnolle \cite[Theorem 8]{Br71})\\
(i)  For Case A, ${\cal C}=\mathbf{R}$ and 0 is a regular point of $\{0\}$. \\
(ii) For Case B, either ${\cal C}=\emptyset$ or ${\cal
C}=\mathbf{R}$, and if
${\cal C}=\mathbf{R}$  then 0 is a regular point of $\{0\}$.\\
 (iii) For Case C, suppose that $X$ is not a compound Poisson process, then\\
 \indent (a) for Case $C_1$, ${\cal C}=\emptyset$;\\
 \indent (b) for Case $C_2$, ${\cal C}=\mathbf{R}^+:=\{x\in \mathbf{R}:x>0\}$ and 0 is not a regular point of $\{0\}$;\\
  \indent(c)  for Case $C_3$, ${\cal C}=\mathbf{R}$ and 0 is not a regular point of $\{0\}$.
 \end{thm}

For Case A, and Case B with ${\cal C}=\mathbf{R}$, only the empty
set is a semipolar set. Hence (H) holds for these two cases. For
Case $C_2$ and Case $C_3$, any singleton $\{x\}$ is semipolar but
non-polar. Thus (H) doesn't hold for these two cases. Therefore,
for one-dimensional L\'evy processes, we need only consider
whether (H) holds for Case B with ${\cal C}=\emptyset$ and Case
$C_1$.

For Case B, Kesten \cite[Theorem 1(f)]{Ke69} tells us that if
$\int_0^{\infty}(1 \wedge x)\mu(dx)<\infty$ or
$\int_{-\infty}^0(1\wedge|x| )\mu(dx)<\infty$, then ${\cal
C}=\mathbf{R}$. Thus, any $x\in \mathbf{R}$ is a regular point of
$\{x\}$ and hence (H) holds for this case. As a consequence, any
spectrally one sided one-dimensional L\'{e}vy process with
unbounded variation satisfies (H). Therefore, for Case B, we need
only consider the case that  both $\int_0^{\infty}(1\wedge x
)\mu(dx)=\infty$ and $\int_{-\infty}^0(1\wedge|x|
)\mu(dx)=\infty$.

Denote by $\mu_{+}$ and $\mu_{-}$ the restriction of the L\'{e}vy
measure $\mu$ on $(0,\infty)$ and $(-\infty,0)$, respectively. Let
$X_1$ and $X_2$ be two independent L\'{e}vy processes with
L\'{e}vy measures $\mu_+$ and $\mu_-$, respectively. For Case B
with $\int_0^{\infty}(1\wedge x )\mu(dx)=\infty$ and
$\int_{-\infty}^0(1\wedge|x| )\mu(dx)=\infty$, both $X_1$ and
$X_2$ belong to Case B with ${\cal C}=\mathbf{R}$ and hence
satisfy (H). Obviously, $X$ can be regarded as the sum of $X_1$
and $X_2$. This observation provides a motivation for us to
consider (H) for the sum of two independent L\'{e}vy processes.

\subsection{Main results}

First, we present a result which extends  \cite[Theorem
1(f)]{Ke69}. Let $\mu$ be the L\'{e}vy measure. We denote by
$\bar{\mu}_{-}$ the image measure of $\mu_{-}$ under the map
$$
x\mapsto -x, \ \forall x\in (-\infty,0).
$$

\begin{thm}\label{thm2.5}
Suppose that $Q=0$ and  $\int_{0}^{\infty}(1\wedge x
)\mu_+(dx)=\infty$. If there exist $\delta\in (0,1), k\in [0,1)$,
and a measure $\nu$ on $\mathbf{R}^+$ satisfying
$\int_{(0,\delta)}x\nu(dx)<\infty$, such that
\begin{eqnarray}\label{pro-a}
\bar{\mu}_{-}\leq k\mu_++\nu.
\end{eqnarray}
Then $X$ satisfies (H).
\end{thm}
{\bf Proof.} We assume without loss of generality that $k>0$.
Define  $\mu_2$ to be the symmetric measure on
$\mathbf{R}\backslash\{0\}$ satisfying
$\mu_2=(\bar{\mu}_{-}-\nu)^+$ on $(0,\delta)$ and $\mu_2=0$ on $[\delta,\infty)$, where
$(\bar{\mu}_{-}-\nu)^+$ denotes the positive part of the signed
measure $\bar{\mu}_{-}-\nu$. Denote $\mu_1=\mu-\mu_2$. Let
$X_1$ and $X_2$ be two independent one-dimensional L\'{e}vy
processes with L\'{e}vy-Khintchine exponents $(a,0,\mu_1)$ and
$(0,0,\mu_2)$, respectively. Since $X$ and $X_1+X_2$ have the same
law, to show that $X$ satisfies (H), it is sufficient to show that
$X_1+X_2$ satisfies (H). We denote by $\psi_1$ and $\psi_2$ the
L\'{e}vy-Khintchine exponents of $X_1$ and $X_2$, respectively.

By  (\ref{pro-a}), we get
$$
\int_0^{\infty}(1\wedge x
)\mu_1(dx)\ge\int_{(0,\delta)}x\mu_1(dx)\geq
(1-k)\int_{(0,\delta)}x\mu_+(dx)=\infty,
$$
and
\begin{eqnarray*}\int_{-\infty}^0(1\wedge|x| )\mu_1(dx)&=&\int_{(-\infty,-\delta]}(1\wedge|x| )\mu_{-}(dx)+\int_{(-\delta,0)}|x|\mu_1(dx)\\
&\leq&\mu_{-}((-\infty,-\delta])+\int_{(0,\delta)}x\nu(dx)\\
&<&\infty.
\end{eqnarray*}
Then, we obtain by \cite[Theorem 1(f)]{Ke69} that $X_1$ belongs to
Case B with ${\cal C}=\mathbf{R}$. Therefore, we obtain by
\cite{Ke69} that
\begin{eqnarray}\label{pro-b}
\int_0^{\infty}{\rm Re}([1+\psi_1(z)]^{-1})dz<\infty.
\end{eqnarray}

By (\ref{pro-a}) and the definition of $\psi_2$, we obtain that
for $z\in \mathbf{R}$,
\begin{eqnarray}\label{pro-c}
\psi_2(z)={\rm Re}\psi_2(z)&=&2\int_{(0,\delta)}(1-\cos (zx))\mu_2(dx)\nonumber\\
&\leq&2k\int_{(0,\delta)}(1-\cos (zx))\mu_+(dx)\nonumber\\
&\le&\frac{2k}{1-k}\int_{(0,\delta)}(1-\cos (zx))\mu_1(dx)\nonumber\\
&\leq &\frac{2k}{1-k}{\rm Re}\psi_1(z).
\end{eqnarray}
By (\ref{pro-b}) and (\ref{pro-c}), we get
\begin{eqnarray*}\label{pro-d}
&&\int_0^{\infty}{\rm Re}([1+\psi_1(z)+\psi_2(z)]^{-1})dz\nonumber\\
&&=\int_0^{\infty}\frac{1}{1+{\rm Re}\psi_1(z)+{\rm Re}\psi_2(z)+\frac{({\rm Im}\psi_1(z))^2}{1+{\rm Re}\psi_1(z)+{\rm Re}\psi_2(z)}}dz\nonumber\\
&&\leq \int_0^{\infty}\frac{1}{1+{\rm Re}\psi_1(z)+\frac{({\rm Im}\psi_1(z))^2}{1+(\frac{1+k}{1-k}){\rm Re}\psi_1(z)}}dz\nonumber\\
&&\leq \frac{1+k}{1-k}\int_0^{\infty}{\rm
Re}([1+\psi_1(z)]^{-1})dz\nonumber\\
&&<\infty.
\end{eqnarray*}
Then, we obtain by \cite{Ke69} that any singleton is non-polar for
$X_1+X_2$. Hence any point $x\in \mathbf{R}$ is a regular point of
$\{x\}$ by Theorem \ref{thm1}(ii). Therefore, $X_1+X_2$ satisfies
(H). \hfill\fbox

We now give a novel condition on the L\'evy measure $\mu$ which implies (H) for a
large class of one-dimensional L\'{e}vy processes.
\begin{thm}\label{thm-2.6}
If
\begin{eqnarray}\label{thm-2.6-a}
\liminf_{\varepsilon\downarrow
0}\frac{\int_{-\varepsilon}^{\varepsilon}x^2\mu(dx)}{\varepsilon/|\log
\varepsilon|}>0,
\end{eqnarray}
then $X$ satisfies (H).
\end{thm}

Note that, different from most sufficient conditions given in the
diagram of Section 1, condition (\ref{thm-2.6-a}) does not require
any controllability of ${\rm Im}(\psi)$ by ${\rm Re}(\psi)$.
Before proving Theorem \ref{thm-2.6}, we give a necessary and
sufficient condition for the validity of (H) for general L\'evy
processes.

\begin{pro}\label{pro123} Suppose that $X$ is a L\'{e}vy
process on $\mathbf{R}^n$  which has resolvent densities w.r.t.
the Lebesgue measure. Let
 $f$ be a positive increasing
function on $[1,\infty)$ such that $ \int_N^{\infty}(\lambda
f(\lambda))^{-1}d\lambda=\infty $ for some $N\ge 1$. Then (H)
holds for $X$ if and only if
\begin{eqnarray*}
\lim_{\lambda\rightarrow\infty}\int_{\{B(z)>A(z)f(A(z))\}}\frac{\lambda}{\lambda^2+B^2(z)}|\hat
\nu(z)|^2dz=0
\end{eqnarray*}
for any finite measure $\nu$ of finite 1-energy.
\end{pro}

\noindent {\bf Proof.} This is a direct consequence of
\cite[Theorems 4.3 and 5.1]{HSZ15}. \hfill\fbox

\vskip 0.3cm

\noindent {\bf Proof of Theorem \ref{thm-2.6}.}
 By (\ref{thm-2.6-a}), we know that there exist
constants $N_1$ and $c$ satisfying $N_1>1$ and $0<c<1$  such that
$$
\int_{-\frac{1}{|z|}}^{\frac{1}{|z|}}x^2\mu(dx)\ge \frac{c}{|z|\log|z|},\ \ {\rm if}\ |z|\geq N_1.
$$
Note that $1-\cos x\geq \frac{x^2}{4}$ when $|x|\leq 1$. Then, for
$|z|\geq N_1$, we have
\begin{eqnarray}\label{thm-2.6-b}
{\rm Re}\psi(z)&=&\int_{\mathbf{R}}(1-\cos(zx))\mu(dx)\nonumber\\
&\geq&\int_{-\frac{1}{|z|}}^{\frac{1}{|z|}}(1-\cos(zx))\mu(dx)\nonumber\\
&\geq&\frac{z^2}{4}\int_{-\frac{1}{|z|}}^{\frac{1}{|z|}}x^2\mu(dx)\nonumber\\
&\geq&\frac{c|z|}{4\log|z|}.
\end{eqnarray}

We define $f(\lambda)=\frac{4}{c}\log(\frac{4
\lambda}{c})[\log\log (\frac{4 \lambda}{c})]$ for $\lambda>c$.
Then, $f(\lambda)$ is a positive increasing function on $(c,\infty)$ and
satisfy
$$
\int_c^{\infty}\frac{1}{\lambda f(\lambda)}d\lambda=\int_{4}^{\infty}\frac{1}{u\log
u[\log\log u]}du=\infty.
$$
We fix a constant $\alpha$ satisfying $0<\alpha<1$. By $\lim_{z\to\infty}\frac{z/\log
z}{z^{\alpha}}=+\infty$, we know that there exists a constant $N_2>0$ such
that
\begin{eqnarray}\label{thm-2.6-c}
\frac{z}{\log z}\geq z^{\alpha},\    \forall z\geq N_2.
\end{eqnarray}
We define $g(z)=\log\log(\frac{z}{\log z})$ for $z>e$. It is easy to see
that $g(z)$ is an increasing positive function on $(e,\infty)$.

By (\ref{thm-2.6-b}) and (\ref{thm-2.6-c}), we obtain that for any
$N_0> \max\{N_1, N_2, e\}$,
\begin{eqnarray*}
&&\limsup_{\lambda\to\infty}\int_{\{B(z)> A(z)f(A(z))\}}\frac{\lambda}{\lambda^2+B^2(z)}dz\\
&&=\limsup_{\lambda\to\infty}\int_{\{B(z)> A(z)f(A(z)),|z|> N_0\}}\frac{\lambda}{\lambda^2+B^2(z)}dz\\
&&\leq \limsup_{\lambda\to\infty}\int_{\{|z|>N_0\}}\frac{\lambda}{\lambda^2+\frac{z^2}{\log^2|z|}\log^2(\frac{|z|}{\log|z|})[\log\log(\frac{|z|}{\log|z|})]^2}dz\\
&&\leq \limsup_{\lambda\to\infty}\int_{\{|z|>N_0\}}\frac{\lambda}{\lambda^2+|z|^2\left(\frac{\log(|z|^{\alpha})}{\log|z|}\right)^2g^2(|z|)}dz\\
&&\leq \lim_{\lambda\to\infty}\int_{\{|z|>N_0\}}\frac{\lambda}{\lambda^2+\alpha^2g^2(N_0)|z|^2}dz\\
&&=\limsup_{\lambda\to\infty}\frac{2}{\alpha g(N_0)}\int_{\alpha N_0g(N_0)}^{\infty}\frac{\lambda}{\lambda^2+u^2}du\\
&&\leq \frac{\pi}{\alpha g(N_0)}.
\end{eqnarray*}
Since  $\lim_{N_0\to\infty}g(N_0)=0$, we obtain
\begin{eqnarray}\label{sdf}
\lim_{\lambda\to\infty}\int_{\{B(z)>
A(z)f(A(z))\}}\frac{\lambda}{\lambda^2+B^2(z)}dz=0.
\end{eqnarray}
By (\ref{thm-2.6-b}) and \cite{HW42}, we know that $X$ has bounded
continuous transition densities. Therefore, $X$ satisfies (H) by
(\ref{sdf}) and Proposition \ref{pro123}.\hfill\fbox

\begin{rem}\label{rem-2.7}
For $\alpha>0$, we define the measure $\nu_{\alpha}$ on $(-1,1)$
by
$$
\nu_{\alpha}(dx):=|x\log |x||^{1+\alpha}\mu(dx),\ x\in (-1,1).
$$
We remark that our condition (\ref{thm-2.6-a}) only requires
slightly more than $\nu_{\alpha}$ is an  infinite measures on
$(-1,1)$ for any $\alpha>0$.

(i) Condition  (\ref{thm-2.6-a})  implies that any $\nu_{\alpha}$
is an infinite measure on $(-1,1)$.  In fact, by
(\ref{thm-2.6-a}), we get
\begin{eqnarray}\label{rem-2.7-b}
\lim_{\varepsilon\downarrow 0}\frac{|\log
\varepsilon|^{1+\alpha}\int_{-\varepsilon}^\varepsilon
x^2\mu(dx)}{\varepsilon}=\infty.
\end{eqnarray}
If  $\nu_{\alpha}$ is a finite measure on $(-1,1)$, then
\begin{eqnarray*}
&&\limsup_{\varepsilon\downarrow 0}\frac{|\log \varepsilon|^{1+\alpha}\int_{-\varepsilon}^\varepsilon x^2\mu(dx)}{\varepsilon}\\
&&\leq \limsup_{\varepsilon\downarrow 0}\frac{\int_{-\varepsilon}^\varepsilon x^2|\log|x||^{1+\alpha}\mu(dx)}{\varepsilon}\\
&&=\limsup_{\varepsilon\downarrow
0}\frac{\int_{-\varepsilon}^\varepsilon|x|\nu_{\alpha}(dx)}{\varepsilon}\\
&&\leq \nu_{\alpha}(-1,1),
\end{eqnarray*}
which contradicts (\ref{rem-2.7-b}).

(ii)  If for some $\beta>2$,
\begin{eqnarray}\label{rem-2.7-c}
\liminf_{\varepsilon\downarrow
0}\frac{\int_{-\varepsilon}^{\varepsilon}x^2\mu(dx)}{\varepsilon/|\log
\varepsilon|^{\beta}}=0,
\end{eqnarray}
then $\nu_{\alpha}$ is a finite measure on $(-1,1)$ for any
$\alpha\in (0,\beta-2)$.

We only prove  $\nu_{\alpha}(0,1)<\infty$. The proof that
$\nu_{\alpha}(-1,0)<\infty$ is similar so we omit it. By
(\ref{rem-2.7-c}), we know that there exist constants $c$ and
$\delta$ satisfying $c>0$ and $0<\delta<1$ such that
$$
\int_0^\varepsilon x^2\mu(dx)\leq \frac{c\varepsilon}{|\log
\varepsilon|^{\beta}},\ \forall \varepsilon\in (0,\delta).
$$

Note that $f(x)=x/|\log x|^{1+\alpha}$ is an increasing function
on $(0,1)$. Then, for any $\varepsilon\in (0,\delta)$, we have
\begin{eqnarray*}
\frac{\varepsilon/2}{|\log
(\varepsilon/2)|^{1+\alpha}}\nu([\varepsilon/2,\varepsilon])\leq
\int_{\frac{\varepsilon}{2}}^\varepsilon\frac{x}{|\log
x|^{1+\alpha}}\nu(dx)\leq \int_0^\varepsilon x^2\mu(dx)\leq
\frac{c\varepsilon}{|\log \varepsilon|^{\beta}},
\end{eqnarray*}
which implies that
$$
\nu([\varepsilon/2,\varepsilon])\leq \frac{2c|\log
(\varepsilon/2)|^{1+\alpha}}{|\log \varepsilon|^{\beta}}.
$$
We fix a $K\in \mathbf{N}$ satisfying $\frac{1}{2^{K}}<\delta$.
Then,
\begin{eqnarray*}
\nu(0,1)&=&\sum_{n=1}^{K}\nu([1/2^n, 1/2^{(n-1)}))+\sum_{n=K+1}^{\infty}\nu([1/2^n, 1/2^{(n-1)}))\\
&\leq &\sum_{n=1}^{K}\nu([1/2^n, 1/2^{(n-1)}))+\sum_{n=K+1}^{\infty}\frac{2c|\log (1/2^n)|^{1+\alpha}}{|\log (1/2^{(n-1)})|^{\beta}}\\
&=&\sum_{n=1}^{K}\nu([1/2^{(n-1)},1/2^n))+2c\sum_{n=K+1}^{\infty}\frac{(n\log
2)^{1+\alpha}}{((n-1)\log 2)^{\beta}}\\ &<&\infty.
\end{eqnarray*}

\end{rem}

From the proof of Theorem \ref{thm-2.6}, we can see that the
following result extending \cite[Theorem 4.7]{HSZ15} holds.
\begin{pro} If
$$
\liminf_{|z|\to\infty}\frac{{\rm Re}\psi(z)}{|z|/\log|z|}>0,
$$
then $X$ satisfies (H).
\end{pro}

Following the proof of Theorem \ref{thm-2.6}, we can also prove
the following proposition.

\begin{pro}
If
\begin{eqnarray*}
\liminf_{\varepsilon\to
0}\frac{\int_{-\varepsilon}^{\varepsilon}x^2\mu(dx)}{\frac{\varepsilon}{|\log
\varepsilon|[\log|\log\varepsilon|]}}>0,
\end{eqnarray*}
 then $X$ satisfies (H).
\end{pro}

\subsection{An example}

We give an application of Theorem \ref{thm-2.6}. Note that in the
following example, there is no assumption on $a$ or $Q$.
\begin{exa}\label{exa-2.9}
Let $X$ be a L\'evy process on $\mathbf{R}$ with L\'evy measure
$\mu$. Suppose that there exist positive constants $c,\delta$, and
a finite measure $\nu$ on $(0,\delta)$ such that
$$
\mu(dx)+\nu(dx)\geq \frac{c}{x^2|\log x|}dx\ \ {\rm on}\
(0,\delta).
$$
Then $X$ satisfies (H).

In fact, we have
\begin{eqnarray*}
\liminf_{\varepsilon\downarrow 0}\frac{\int_0^\varepsilon
x^2\cdot\frac{c}{x^2|\log x|}dx}{\varepsilon/|\log \varepsilon|}
&\geq& \liminf_{\varepsilon\downarrow 0}\frac{\int_{\varepsilon/2}^\varepsilon \frac{c}{|\log x|}dx}{\varepsilon/|\log \varepsilon|}\nonumber\\
&\geq& \liminf_{\varepsilon\downarrow 0}\frac{ \frac{c}{|\log \varepsilon|}\cdot \frac{\varepsilon}{2}}{\varepsilon/|\log \varepsilon|}\nonumber\\
&=&\frac{c}{2},
\end{eqnarray*}
and
\begin{eqnarray*}
\limsup_{\varepsilon\downarrow 0}\frac{\int_0^\varepsilon
x^2\nu(dx)}{\varepsilon/|\log \varepsilon|} &\leq
&\limsup_{\varepsilon\downarrow
0}\frac{\varepsilon^2\nu(0,1)}{\varepsilon/|\log \varepsilon|}=0.
\end{eqnarray*}
Then (\ref{thm-2.6-a}) holds and therefore $X$ satisfies (H) by Theorem
\ref{thm-2.6}.
\end{exa}

%
%
%

\section{(H) for sum of L\'{e}vy processes: no assumption on resolvent densities}\setcounter{equation}{0}

From now on till the end of the paper, we consider Hunt's
hypothesis (H) for general $\mathbf{R}^n$-valued L\'{e}vy
processes. In this section, we discuss (H) for the sum of two
independent L\'{e}vy processes without any assumption on resolvent
densities. In the next section, we discuss (H) for the sum of two
independent L\'{e}vy processes under the assumption that resolvent
densities exist.
\subsection{Main results}

\begin{thm}\label{thm-3.1}
Let $X_1$ and $X_2$ be two independent L\'evy processes on
$\mathbf{R}^n$. If $X_1$ satisfies (H) and $X_2$ is a compound
Poisson process, then $X_1+X_2$ satisfies (H).
\end{thm}

\begin{thm}\label{thm-3.2}
Let $X_1$ and $X_2$ be  two independent L\'evy processes on
$\mathbf{R}^n$. If both $X_1$ and $X_2$ satisfy condition (S),
then $X_1+X_2$ satisfies (H).
\end{thm}

As a direct consequence of Theorem  \ref{thm-3.1}, we can
strengthen \cite[Theorem 2.1]{HSZ15} as follows:

\begin{pro}\label{pro-3.5} Let $X$ be a L\'{e}vy process  on $\mathbf{R}^n$ with
L\'{e}vy-Khintchine exponent  $(a,Q,\mu)$. Suppose  that $\mu_1$
is a finite measure on $\mathbf{R}^n\backslash\{0\}$ such that
$\mu_1\leq \mu$. Denote $\mu_2:=\mu-\mu_1$ and let $X'$ be a
L\'{e}vy process on $\mathbf{R}^n$ with L\'{e}vy-Khintchine
exponent $(a',Q,\mu_2)$, where $a':=a+\int_{\{|x|<1\}}x\mu_1(dx).$
Then

(i) $X$ and $X'$ have same semipolar sets.

(ii) $X$ and $X'$ have same essentially polar sets.

(iii)  if $X'$ satisfies (H), then $X$ satisfies (H).

(iv) if $X$ satisfies (H) and $X'$ has resolvent densities w.r.t. the Lebesgue measure, then
$X'$ satisfies (H).
\end{pro}

\subsection{Proof of Theorem \ref{thm-3.1}}

Before proving Theorem \ref{thm-3.1}, we present some lemmas,
which have their own interests.

\begin{lem}\label{lem-3.2}
Let $X$ be a L\'{e}vy process on $\mathbf{R}^n\ (n>1)$ satisfying
(H). Then, for any nonempty proper subspace $S$ of $\mathbf{R}^n$,
the projection process $Y$ of $X$ on $S$ satisfies (H).
 \end{lem}
{\bf Proof.} By virtue of the orthogonal transformation (cf.
\cite[Section 2.2]{HS12}),  we can assume without loss of
generality that $S=\{(x_1,\cdots,x_n)\in
\mathbf{R}^n|x_{k+1}=\cdots=x_n=0\}$ for some integer $k,1\leq
k<n$. Then, the projection process $Y$ of $X$ can be regarded as a
L\'{e}vy process on $\mathbf{R}^k$. Let $C\subset \mathbf{R}^k$ be
a semipololar set for $Y$. We define
$$
D=\{(x_1,\cdots,x_n)\in \mathbf{R}^n|(x_1,\cdots,x_k)\in C\}.
$$
By the definition of semipolar set, we find that $D$ is a
semipolar set for $X$. Further, by the assumption that $X$
satisfies (H), we conclude that $D$ is a  polar set for $X$.
Therefore, as the projection of $D$ on $S$, $C$ is a polar set for
$Y$. \hfill\fbox

\begin{lem}\label{lem-3.3}
Let $X$ be a L\'{e}vy process on $\mathbf{R}^n\ (n>1)$  with
L\'{e}vy-Khintchine exponent $(a,Q,\mu)$. Suppose that for some
proper subspace $S$ of $\mathbf{R}^n$, the projection process
$X_S$ of $X$ on $S$ satisfies (H) and $\mu(\mathbf{R}^n\backslash
S)<\infty$.  Then $X$ satisfies (H).
\end{lem}

\noindent{\bf Proof.} By virtue of the orthogonal transformation,
we can assume without loss of generality that
$S=\{(x_1,\cdots,x_n)\in \mathbf{R}^n|x_{k+1}=\cdots=x_n=0\}$ for
some integer $k$, $1\leq k<n$.  By the L\'{e}vy-It\^{o} decomposition
(cf. the proof of \cite[Theorem 1.2]{HS12}), we may express $X$ as
$$
X=X^{(1)}+X^{(2)},
$$
where $X^{(1)}=(X_S,0)$ can be regarded as a $k$-dimensional
L\'{e}vy process on $\mathbf{R}^k\times \{0\}$ which satisfies
(H), and $X^{(2)}$ is a compound Poisson process on $\mathbf{R}^n$
which is independent of $X^{(1)}$. Then, by following the proof of
(ii) $\Rightarrow$ (i) of \cite[Theorem 1.2]{HS12}, we conclude
that $X$ satisfies (H).\hfill\fbox

\begin{lem}\label{lem-3.4}
Let $X_1$ and $X_2$ be two independent L\'{e}vy processes on
$\mathbf{R}^m$ and $\mathbf{R}^n$, respectively. If $X_1$
satisfies (H) and $X_2$ is a compound Poisson process, then
$X=(X_1,X_2)$ satisfies (H).
\end{lem}
{\bf Proof.} This is a direct consequence of Lemma
\ref{lem-3.3}.\hfill\fbox

\bigskip

\noindent{\bf Proof of Theorem \ref{thm-3.1}}. By Lemma
\ref{lem-3.4}, we find that the $\mathbf{R}^{2n}$-valued L\'{e}vy
process $(X_1,X_2)$ satisfies (H). Further, by the orthogonal
transformation, we find that the L\'evy process
$\frac{1}{\sqrt{2}}(X_1+X_2,X_2-X_1)$ satisfies (H). Therefore,
$X_1+X_2$ satisfies (H) by Lemma \ref{lem-3.2}.\hfill\fbox

\subsection{Proof of Theorem \ref{thm-3.2}}

Before giving the proof for Theorem \ref{thm-3.2}, we prove the
following lemma.

\begin{lem}\label{lem-3.7}
 Let $M$ be a symmetric nonnegative definite $n\times n$ matrix. Then, $x\in \sqrt{M}\mathbf{R}^n$ if and only if there exists a constant $c>0$ such that
 \begin{eqnarray}\label{lem1-a}
 |\langle x,z\rangle|\leq c\sqrt{\langle z,Mz\rangle},\ \forall z\in \mathbf{R}^n.
 \end{eqnarray}
 \end{lem}
 \noindent{\bf Proof.} Suppose that $x\in \sqrt{M}\mathbf{R}^n$. Then, there exists a $y\in \mathbf{R}^n$ such that $x=\sqrt{M}y$ and thus
 \begin{eqnarray*}
 |\langle x,z\rangle|&=&|\langle \sqrt{M}y,z\rangle|\\
 &=&|\langle y,\sqrt{M}z\rangle|\\
 &\leq&\sqrt{\langle y,y\rangle}\sqrt{\langle\sqrt{M}z,\sqrt{M}z\rangle}\\
 &=&\sqrt{\langle y,y\rangle}\sqrt{\langle z,Mz\rangle}.
 \end{eqnarray*}
 Therefore, (\ref{lem1-a}) holds with $c=1+\sqrt{\langle y,y\rangle}$.

Now we suppose that (\ref{lem1-a}) holds.  Denote by $k$ the rank
of $M$. If $k=n$ or 0, it is easy to see that $x\in
\sqrt{M}\mathbf{R}^n$. Hence we may assume that $n\ge 2$ and
$1\leq k\le n-1$.
 Since $M$ is a symmetric nonnegative definite $n\times n$ matrix, there exists an orthogonal matrix $O$ such that
$$
OMO^T=diag(\lambda_1,\dots,\lambda_n):=F,
$$
where $\lambda_1\geq  \cdots\geq \lambda_k> 0,\lambda_i=0$ for
$i=k+1,\cdots,n$, and $O^T$ denotes the transpose of $O$. We can
rewrite (\ref{lem1-a}) as follows:
\begin{eqnarray*}
|\langle Ox,Oz\rangle|\leq c\sqrt{\langle Oz,F(Oz)\rangle},\
\forall z\in \mathbf{R}^n,
\end{eqnarray*}
equivalently,
\begin{eqnarray}\label{lem1-b}
|\langle Ox,z'\rangle|\leq c\sqrt{\langle z',Fz'\rangle},\ \forall
z'\in \mathbf{R}^n.
\end{eqnarray}

We claim that $Ox\in \sqrt{F}\mathbf{R}^n=\mathbf{R}^k\times
\{0\}$. Let $Ox=(\bar{x}_1,\bar{x}_2,\cdots,\bar{x}_n)$. If
$Ox\notin \mathbf{R}^k\times\{0\}$, then there exists $j\in
\{k+1,\cdots,n\}$ such that $\bar{x}_j\neq 0$. Let
$z'=(z'_1,\cdots,z'_n)$ with $z'_j=1$ and $z'_i=0$ for $i\not=0$.
Thus, we obtain  by (\ref{lem1-b})  that
$$
0<|\bar{x}_j|=|\langle Ox,z'\rangle|\leq c\sqrt{\langle
z',Fz'\rangle}=0.
$$
This is a contradiction and hence $Ox\in \sqrt{F}\mathbf{R}^n$.
Therefore, $x\in \sqrt{M}\mathbf{R}^n$. \hfill\fbox
\bigskip

\noindent{\bf Proof of Theorem \ref{thm-3.2}.} We denote the
L\'{e}vy-Khintchine exponents of $X_1$ and $X_2$ by
$(a_1,Q_1,\mu_1)$ and $(a_2,Q_2,\mu_2)$, respectively. By Lemma
\ref{lem-3.7}, we find that $\sqrt{Q_1}\mathbf{R}^n\subset
\sqrt{Q_1+Q_2}\mathbf{R}^n$ and $\sqrt{Q_2}\mathbf{R}^n\subset
\sqrt{Q_1+Q_2}\mathbf{R}^n$. Thus
\begin{equation}\label{hjk}(\mu_1+\mu_2)(\mathbf{R}^n\backslash
\sqrt{Q_1+Q_2}\mathbf{R}^n)<\infty.\end{equation}

By \cite[Theorem 1.2]{HS12}, we know that both $X_1$ and $X_2$
satisfy the Kanda-Forst condition and hence $X_1+X_2$ satisfies
the Kanda-Forst condition. Therefore,  $X_1+X_2$ satisfies (H) by
(\ref{hjk}) and \cite[Theorem 1.2]{HS12}.\hfill\fbox

\section{(H) for sum of L\'{e}vy processes under assumption that resolvent densities exist}\setcounter{equation}{0}
Throughout this section, we assume that $X_1$ and $X_2$ are two
independent L\'{e}vy processes on $\mathbf{R}^n$ such that
$X_1+X_2$ has resolvent densities w.r.t. the Lebesgue measure. We
denote by $\psi_1$ and $\psi_2$ the L\'{e}vy-Khintchine exponents
of $X_1$ and $X_2$, respectively.

\subsection{Main results}

\begin{thm}\label{thm-4.1}
Suppose that

(i) $X_1$ has resolvent densities w.r.t. the Lebesgue measure and
satisfies (H).

(ii) Any finite measure $\nu$ of finite 1-energy w.r.t. $X_1+X_2$
has finite 1-energy w.r.t. $X_1$.

(iii) There exists a constant $c>0$ such that
\begin{eqnarray*}
|{\rm Im}\psi_2|\leq c(1+{\rm Re}\psi_1+{\rm Re}\psi_2).
\end{eqnarray*}
Then $X_1+X_2$ satisfies (H).
\end{thm}
\begin{pro}\label{pro-3.12} If one of the following conditions is fulfilled, then any finite measure $\nu$ of finite 1-energy w.r.t. $X_1+X_2$
has finite 1-energy w.r.t. $X_1$.

(i) There exists a constant $c>0$ such that
$$
|\psi_2|\leq c(1+{\rm Re}(\psi_1)).
$$

(ii) There exists a constant $c>0$ such that
\begin{eqnarray*}
\left\{
\begin{array}{l}
{\rm Re}\psi_2\leq c\left(1+{\rm Re}\psi_1+\frac{({\rm Im}\psi_1)^2}{1+{\rm Re}\psi_1}\right),\\
|{\rm Im}\psi_2|\leq c(1+{\rm Re}\psi_1+{\rm Re}\psi_2).
\end{array}\right.
\end{eqnarray*}

(iii) There exists a constant $c>0$ such that
\begin{eqnarray}\label{pro-3.12-a}
\left\{
\begin{array}{l}
{\rm Re}\psi_2\leq c\left(1+{\rm Re}\psi_1+\frac{({\rm Im}\psi_1)^2}{1+{\rm Re}\psi_1)}\right),\\
({\rm Im}\psi_2)^2\leq c(1+{\rm Re}\psi_1+{\rm
Re}\psi_2)\left(1+{\rm Re}\psi_1+\frac{({\rm Im}\psi_1)^2}{1+{\rm
Re}\psi_1}\right).
\end{array}\right.
\end{eqnarray}
\end{pro}

\begin{cor}\label{cor-4.2}
Suppose that

(i) $X_1$ has bounded resolvent densities w.r.t. the Lebesgue
measure and satisfies (H).

(ii)  There exists a constant $c>0$ such that
$$
|{\rm Im}\psi_2|\leq c(1+{\rm Re}\psi_1+{\rm Re}\psi_2).
$$
Then $X_1+X_2$ satisfies (H).
\end{cor}

\begin{rem}\label{rem4.10}
Let $X$ be a one-dimensional L\'{e}vy process and  the set ${\cal
C}$ be defined as in (\ref{sec2-a}).  By \cite[Theorem 43.21, Case
5]{Sa99}, we know that if $X$ belongs to  {\it Case B} (defined as
in Section 2) with ${\cal C}=\mathbf{R}$, then $X$ has bounded
resolvent densities w.r.t the Lebesgue measure.  In particular,

 (i) the one-dimensional Brownian motion has bounded resolvent
 densities.

 (ii) any spectrally one sided one-dimensional L\'{e}vy process with unbounded variation has bounded resolvent
 densities.

(iii)  any one-dimensional L\'{e}vy process satisfying the
conditions of Theorem \ref{thm2.5} has bounded resolvent
 densities.
\end{rem}

\begin{pro} \label{pro-4.3}
Let $f$ be a positive increasing function on $[1,\infty)$ such
that $ \int_N^{\infty}(\lambda f(\lambda))^{-1}d\lambda=\infty $
for some $N\ge 1$.  Suppose that

(i) There are two measurable functions $\phi_{11}$ and $\phi_{12}$
on $\mathbf{R}^n$ such that $\rm{Im}\psi_1=\phi_{11}+\phi_{12}$,
and
$$
|\phi_{11}|\leq (1+{\rm Re}\psi_1)f(1+{\rm Re}\psi_1),\ \
\int_{{\mathbf{R}^n}}\frac{|\phi_{12}(z)|}{|1+\psi_1(z)|^2}dz<\infty.
$$

(ii)
\begin{eqnarray*}\label{pro2-b}
|{\rm Im}\psi_2|\leq (1+{\rm Re}\psi_1+{\rm Re}\psi_2)f(1+{\rm
Re}\psi_1+{\rm Re}\psi_2).
\end{eqnarray*}
Then $X_1+X_2$ satisfies (H).
 \end{pro}

\subsection{Proofs}

Before giving the proof for Theorem \ref{thm-4.1}, we prove the
following lemma.
\begin{lem}\label{lem-3.14}
Suppose that there exists a constant $c>0$ such that
\begin{eqnarray}\label{lem2.7-a}
|{\rm Im}\psi_2|\leq c(1+{\rm Re}\psi_1+{\rm Re}\psi_2).
\end{eqnarray}
Then, there exists a constant $\gamma>0$ such that
\begin{eqnarray*}
|1+\psi_1+\psi_2|^2\geq \gamma |1+\psi_1|^2.
\end{eqnarray*}
\end{lem}
{\bf Proof.} Suppose that (\ref{lem2.7-a}) holds. We take
$\gamma\in (0,\frac{1}{4})$ such that
$(1-\gamma)(1+\frac{1}{4c^2})>1$. Then, for any $x\in \mathbf{R}$,
we have
\begin{eqnarray*}
(x+1)^2-\left(\gamma x^2-\frac{1}{4c^2}\right)&=&(1-\gamma)x^2+2x+\left(1+\frac{1}{4c^2}\right)\nonumber\\
&=&(1-\gamma)\left(x+\frac{1}{1-\gamma}\right)^2+\frac{1}{1-\gamma}
\left((1-\gamma)\left(1+\frac{1}{4c^2}\right)-1\right)\nonumber\\
&\geq&\frac{1}{1-\gamma}\left((1-\gamma)\left(1+\frac{1}{4c^2}\right)-1\right)\\
&>&0,
\end{eqnarray*}
which implies that
\begin{eqnarray}\label{lem2.7-b}
(x+1)^2>\gamma x^2-\frac{1}{4c^2},\ \forall x\in \mathbf{R}.
\end{eqnarray}

By (\ref{lem2.7-b}), we get
\begin{eqnarray}\label{lem2.7-c}
({\rm Im}\psi_1+{\rm Im}\psi_2)^2\geq \gamma ({\rm
Im}\psi_1)^2-\frac{1}{4c^2}({\rm Im}\psi_2)^2.
\end{eqnarray}
Therefore, we obtain by (\ref{lem2.7-a}) and (\ref{lem2.7-c}) that
\begin{eqnarray*}\label{thm4.7-c}
|1+\psi_1+\psi_2|^2&=&(1+{\rm Re}\psi_1+{\rm Re}\psi_2)^2+({\rm Im}\psi_1+{\rm Im}\psi_2)^2\nonumber\\
&=&\left[\left(\frac{1}{2}+\frac{1}{2}{\rm Re}\psi_1\right)+\left(\frac{1}{2}+\frac{1}{2}{\rm Re}\psi_1+{\rm Re}\psi_2\right)\right]^2+({\rm Im}\psi_1+{\rm Im}\psi_2)^2\nonumber\\
&\geq&\left(\frac{1}{2}+\frac{1}{2}{\rm Re}\psi_1\right)^2+\left(\frac{1}{2}+\frac{1}{2}{\rm Re}\psi_1+{\rm Re}\psi_2\right)^2+\gamma ({\rm Im}\psi_1)^2-\frac{1}{4c^2}({\rm Im}\psi_2)^2\quad\quad\nonumber\\
&\geq&\frac{1}{4}\left(1+{\rm Re}\psi_1\right)^2+\frac{1}{4}\left(1+{\rm Re}\psi_1+{\rm Re}\psi_2\right)^2+\gamma ({\rm Im}\psi_1)^2\nonumber\\
&&-\frac{1}{4c^2}\cdot c^2(1+{\rm Re}\psi_1+{\rm Re}\psi_2)^2\nonumber\\
&=&\frac{1}{4}\left(1+{\rm Re}\psi_1\right)^2+\gamma ({\rm Im}\psi_1)^2\nonumber\\
&>&\gamma\left[\left(1+{\rm Re}\psi_1\right)^2+({\rm
Im}\psi_1)^2\right]\\
&=&\gamma |1+\psi_1|^2.
\end{eqnarray*}
The proof is complete.\hfill\fbox

\bigskip \noindent{\bf Proof of Theorem \ref{thm-4.1}}.  Let $\nu$ be a finite measure of finite 1-energy w.r.t. $X_1+X_2$. By Assumption (ii), $\nu$ has finite 1-energy w.r.t. $X_1$.
Then, by Assumption (i) and \cite[Proposition 2.2]{HS16}, we get
\begin{equation}\label{add0}
\lim_{\lambda\rightarrow\infty}
\int_{\mathbf{R}^n}\frac{\lambda}{\lambda^2+|1+\psi_1(z)|^2}|\hat{\nu}(z)|^2dz=0.
\end{equation}
By Assumption (iii) and Lemma \ref{lem-3.14}, we find that there
exists a constant $\gamma>0$ such that
\begin{equation}\label{add}
|1+\psi_1+\psi_2|^2\geq \gamma |1+\psi_1|^2.
\end{equation}
By (\ref{add0}) and (\ref{add}), we obtain that
\begin{eqnarray*} &&\limsup_{\lambda\rightarrow\infty}
\int_{\mathbf{R}^n}\frac{\lambda}{\lambda^2+|1+\psi_1(z)+\psi_2(z)|^2}|\hat{\nu}(z)|^2dz\\
&&\leq \limsup_{\lambda\rightarrow\infty}
\int_{\mathbf{R}^n}\frac{\lambda}{\lambda^2+\gamma|1+\psi_1(z)|^2}|\hat{\nu}(z)|^2dz\\
&&=\limsup_{\lambda\rightarrow\infty}\frac{1}{\sqrt{\gamma}}
\int_{\mathbf{R}^n}\frac{\frac{1}{\sqrt{\gamma}}\lambda}{(\frac{1}{\sqrt{\gamma}}\lambda)^2
+ |1+\psi_1(z)|^2}|\hat{\nu}(z)|^2dz\\
&&=0.
\end{eqnarray*}
Therefore, $X_1+X_2$ satisfies (H) by \cite[Proposition
2.2]{HS16}.\hfill\fbox

\bigskip

\noindent{\bf Proof of Proposition \ref{pro-3.12}.}  It is easy to
see that condition (i) $\Rightarrow$ condition (ii) $\Rightarrow$
condition (iii). In the following, we will prove that if condtion
(iii) is fulfilled, then any finite measure $\nu$ of finite
1-energy w.r.t. $X_1+X_2$ has finite 1-energy w.r.t. $X_1$.

We denote by $\psi$ the L\'{e}vy-Khintchine exponent of $X_1+X_2$.
Suppose that $\nu$ is a finite measure of finite 1-energy w.r.t.
$X_1+X_2$, i.e.,
\begin{eqnarray}\label{pro-3.12-b}
\int_{\mathbf{R}^n}\frac{1+{\rm
Re}\psi(z)}{|1+\psi(z)|^2}|\hat{\nu}(z)|^2dz=\int_{\mathbf{R}^n}{\rm
Re}\left(\frac{1}{1+\psi(z)}\right)|\hat{\nu}(z)|^2dz<\infty.\ \
\end{eqnarray}
By (\ref{pro-3.12-a}), for any $z\in \mathbf{R}^n$, we have
 \begin{eqnarray}\label{pro-3.12-c}
 &&{\rm Re}\left(\frac{1}{1+\psi(z)}\right)=\frac{1}{1+{\rm Re}\psi_1(z)+{\rm Re}\psi_2(z)+\frac{({\rm Im}\psi_1(z)+{\rm Im}\psi_2(z))^2}{1+{\rm Re}\psi_1(z)+{\rm Re}\psi_2(z)}}\nonumber\\
 &&\geq \frac{1}{1+{\rm Re}\psi_1(z)+{\rm Re}\psi_2(z)+\frac{2({\rm Im}\psi_1(z))^2+2({\rm Im}\psi_2(z))^2}{1+{\rm Re}\psi_1(z)+{\rm Re}\psi_2(z)}}\nonumber\\
 &&\geq\frac{1}{1+{\rm Re}\psi_1(z)+\frac{2({\rm Im}\psi_1(z))^2}{1+{\rm Re}\psi_1(z)}
 +c\left(1+{\rm Re}\psi_1(z)+\frac{({\rm Im}\psi_1(z))^2}{1+{\rm Re}\psi_1(z)}\right)
 +\frac{2c(1+{\rm Re}\psi_1(z)+{\rm Re}\psi_2(z))\left(1+{\rm Re}\psi_1(z)+\frac{({\rm Im}\psi_1(z))^2}{1+{\rm Re}\psi_1(z)}\right)}{1+{\rm Re}\psi_1(z)+{\rm Re}\psi_2(z)}}\nonumber\\
   &&= \frac{1}{(1+3c)+(1+3c){\rm Re}\psi_1(z)+(2+3c)\frac{({\rm Im}\psi_1(z))^2}{1+{\rm Re}\psi_1(z)}}\nonumber\\
  &&\geq \frac{1}{2+3c}\cdot \frac{1}{1+{\rm Re}\psi_1(z)+\frac{({\rm Im}\psi_1(z))^2}{1+{\rm Re}\psi_1(z)}}\nonumber\\
  &&=\frac{1}{2+3c}{\rm Re}\left(\frac{1}{1+\psi_1(z)}\right).
 \end{eqnarray}
By (\ref{pro-3.12-b}) and (\ref{pro-3.12-c}), we obtain that
\begin{eqnarray*}\label{pro3.4-b}
\int_{\mathbf{R}^n}{\rm
Re}\left(\frac{1}{1+\psi_1(z)}\right)|\hat{\nu}(z)|^2dz
&\leq&(2+3c)\int_{\mathbf{R}^n}{\rm
Re}\left(\frac{1}{1+\psi(z)}\right)|\hat{\nu}(z)|^2dz<\infty.
\end{eqnarray*}
Therefore, $\nu$ has finite 1-energy w.r.t. $X_1$.\hfill\fbox

\bigskip

\noindent{\bf Proof of Corollary \ref{cor-4.2}}. We denote by
$U^1_{X_1}$ the 1-resolvent of $X_1$. By Assumption (i), for any
finite measure $\nu$, $U_{X_1}^1\nu$ is bounded. Hence
$U_{X_1}^1\nu$ has finite 1-energy w.r.t. $X_1$ by
\cite[Remark]{R88}. The corollary is therefore a direct
consequence of Theorem \ref{thm-4.1}.

\bigskip

\noindent{\bf Proof of Proposition \ref{pro-4.3}.}
We define $A(z)=1+{\rm Re}\psi(z)$ and $B(z)=|1+\psi(z)|$ for
$z\in {\mathbf{R}^n}$. Then $A(z)=1+{\rm Re}\psi_1(z)+{\rm
Re}\psi_2(z)$ and $B(z)=|1+\psi_1(z)+\psi_2(z)|$.
We assume without loss of generality that $f(1)=1/3$. Note that
$B(z)>3\sqrt{2}A(z)f(A(z))$ implies that $|{\rm Im}\psi(z)|> A(z)$
and $|{\rm Im}\psi(z)|> B(z)/\sqrt{2}$. Since $|{\rm
Im}\psi_2|\leq A(z)f(A(z))$, we know that if $|{\rm
Im}\psi(z)|>3A(z)f(A(z))$, then $|{\rm Im}\psi_1(z)|>2A(z)f(A(z))$
and hence $|{\rm Im}\psi_1(z)|\geq 2|{\rm Im}\psi_2(z)|$. Thus
\begin{eqnarray*}
({\rm Im}\psi(z))^2&=&({\rm Im}\psi_1(z)+{\rm Im}\psi_2(z))^2\\
&\geq&(|{\rm Im}\psi_1(z)|-|{\rm Im}\psi_2(z)|)^2\\
&\geq& \frac{1}{4}({\rm Im}\psi_1(z))^2.
\end{eqnarray*}
Note that $|{\rm Im}\psi_1(z)|>2A(z)f(A(z))$ implies that $|{\rm
Im}\psi_1(z)|>\frac{2}{3}(1+{\rm Re}\psi_1(z))$ and  $|\phi_{12}(z)|\geq |{\rm
Im}\psi_1(z)|/2$. Then,  by  the fact that $A(z)\leq c(1+|z|^2)$
for some constant $c>0$ and  the dominated convergence theorem, we
obtain that
\begin{eqnarray*}
&&\sum_{k=1}^{\infty}\int_{\{B(z)>3\sqrt{2}A(z)f(A(z)),\,k\le \frac{|{\rm Im}\psi(z)|}{A(z)}<k+1,\,A(z)\le\lambda<(k+1)|{\rm Im}\psi(z)|\}}\frac{\lambda}{\lambda^2+({\rm Im}\psi(z))^2}|\hat \nu(z)|^2dz\\
&\le&\sum_{k=1}^{\infty}\int_{\{|{\rm Im}\psi(z)|>3A(z)f(A(z)),\,k\le \frac{|{\rm Im}\psi(z)|}{A(z)}<k+1,\,A(z)\le\lambda<(k+1)|{\rm Im}\psi(z)|\}}\frac{\lambda}{\lambda^2+\frac{1}{4}({\rm Im}\psi_1(z))^2}|\hat \nu(z)|^2dz\\
&\le&\sum_{k=1}^{\infty}\int_{\{|{\rm Im}\psi(z)|>3A(z)f(A(z)),\,k\le \frac{|{\rm Im}\psi(z)|}{A(z)}<k+1,\,A(z)\le\lambda<(k+1)|{\rm Im}\psi(z)|\}}\frac{2}{|{\rm Im}\psi_1(z)|}|\hat \nu(z)|^2dz\\
&\le&\sum_{k=1}^{\infty}\int_{\{|{\rm Im}\psi_1(z)|>2A(z)f(A(z)),\,k\le \frac{|{\rm Im}\psi(z)|}{A(z)}<k+1,\,A(z)\le\lambda<(k+1)|{\rm Im}\psi(z)|\}}\frac{4|{\rm Im}\psi_1(z)|}{2({\rm Im}\psi_1(z))^2}| \hat \nu(z)|^2dz\\
&\le&\sum_{k=1}^{\infty}\int_{\{k\le \frac{|{\rm Im}\psi(z)|}{A(z)}<k+1,\,A(z)\le\lambda<(k+1)|{\rm Im}\psi(z)|\}}\frac{8|\phi_{12}(z)|}{(\frac{2}{3}(1+{\rm Re}\psi_1(z))^2+({\rm Im}\psi_1(z))^2}|\hat \nu(z)|^2dz\\
&\le&\sum_{k=1}^{\infty}\int_{\{k\le \frac{|{\rm Im}\psi(z)|}{A(z)}<k+1,\,\lambda<(k+1)^2A(z)\}}\frac{18|\phi_{12}(z)|}{(1+{\rm Re}\psi_1(z))^2+({\rm Im}\psi_1(z))^2}|\hat \nu(z)|^2dz\\
&\le&\sum_{k=1}^{\infty}\int_{\{k\le \frac{|{\rm Im}\psi(z)|}{A(z)}<k+1,\,\lambda<c(k+1)^2(1+|z|^2)\}}\frac{18|\phi_{12}(z)|}{(1+{\rm Re}\psi_1(z))^2+({\rm Im}\psi_1(z))^2}|\hat \nu(z)|^2dz\\
&\rightarrow&0\ \ {\rm as}\ \lambda\rightarrow\infty.
\end{eqnarray*}
Therefore, $X_1+X_2$ satisfies (H) by \cite[Theorem 4.3]{HSZ15}.
\hfill\fbox

\subsection{Examples}

\begin{exa}\label{exa-5.10}
Let  $X_1$ and $X_2$ be two independent L\'{e}vy processes on
$\mathbf{R}^n$. We denote by $\psi_1$ and $\psi_2$ the
L\'{e}vy-Khintchine exponents of $X_1$ and $X_2$, respectively.
Following Blumenthal and Getoor \cite{BG61}, we define the
indices:
\begin{eqnarray*}
\beta_1^{''}&:=&\sup\left\{\alpha\geq 0: \frac{{\rm Re}\psi_1(z)}{|z|^{\alpha}}\to\infty\ \mbox{\rm as}\ |z|\to\infty\right\},\nonumber\\
\beta_2&:=&\inf\left\{\alpha>0: \int_{\{|x|<
1\}}|x|^{\alpha}\nu_2(dx)<\infty\right\},
\end{eqnarray*}
where $\nu_2$ is the L\'{e}vy measure of $X_2$. We will prove
below that if $X_1$ satisfies (H) and $\beta_2<\beta_1^{''}$, then
$X_1+X_2$ satisfies (H).

We fix a $\beta\in (\beta_2,\beta_1^{''})$. Then
\begin{eqnarray}\label{exa-5.10-b} \lim_{|z|\to\infty}\frac{{\rm
Re}\psi_1(z)}{|z|^{\beta}}=\infty.
\end{eqnarray}By \cite[Theorem
3.2]{BG61}, we get
\begin{eqnarray}\label{exa-5.10-a}
\lim_{|z|\to\infty}\frac{|\psi_2(z)|}{|z|^{\beta}}=0.
\end{eqnarray}
(\ref{exa-5.10-b}) and (\ref{exa-5.10-a}) imply  that there
exists a constant $c>0$ such that
$$
|\psi_2(z)|\leq c(1+{\rm Re}\psi_1(z)),\ \forall z\in
\mathbf{R}^n.
$$
By the assumption that $\beta_2<\beta_1^{''}$, we get
$\beta_1^{''}>0$. By (\ref{exa-5.10-b}) and \cite{HW42}, we know that $X_1$ and hence
$X_1+X_2$ have transition densities. Therefore, $X_1+X_2$
satisfies (H) by Theorem \ref{thm-4.1} and Proposition
\ref{pro-3.12}.
\end{exa}

\begin{exa}
Suppose that $\mu$ is a L\'evy measure on $(0,\infty)$ satisfying
$\int_{(0,1)}x\mu(dx)=+\infty$, $\nu$ is a symmetric L\'{e}vy
measure on $\mathbf{R}\backslash \{0\}$, and $a\in \mathbf{R}$.
Let $X$ be a L\'{e}vy process on $\mathbf{R}$ with the
L\'{e}vy-Khintchine exponent $(a,0,\mu+\nu)$.

(i) If $\int_{\{|x|<1\}}|x|\nu(dx)<\infty$, then $X$ satisfies (H)
by Kesten \cite[Theorem 1(f)]{Ke69}.

(ii) If $\int_{\{|x|< 1\}}|x|\nu(dx)=\infty$ and  the restriction
of $\mu$ on $(0,\delta)$ is absolutely continuous w.r.t. the
Lebesgue measure for some constant $\delta$ ($0<\delta<1$), then
$X$ satisfies (H). In fact, let $X_1$ be a L\'{e}vy process on
$\mathbf{R}$ with the L\'{e}vy-Khintchine exponent $(a,0,\mu)$.
Then, $X_1$ has transition densities (cf. \cite[Theorem
27.7]{Sa99}) and bounded resolvent densities (see Remark
\ref{rem4.10}(ii)). It follows that $X$ has transition densities.
Therefore, $X$ satisfies (H) by Corollary \ref{cor-4.2}.
\end{exa}

Before presenting the next example, we recall the definition of type-$(\alpha,\beta)$ subordinator which is introduced in \cite{HS16}.

\begin{defi} (\cite[Definition 4.1]{HS16})
Let $0<\alpha<\beta<1$. A pure jump subordinator $X$ is said to be
of type-$(\alpha,\beta)$ if the L\'evy measure of $X$ has density,
which is denoted by $\rho$, and there exists a constant $c>1$ such
that
\begin{eqnarray*}
\frac{1}{cx^{1+\alpha}}\le\rho(x)\le \frac{c}{x^{1+\beta}},\ \
\forall x\in (0,1].
\end{eqnarray*}
\end{defi}

Up to now it is still unknown if any pure jump subordinator of type-$(\alpha,\beta)$  satisfies (H). In \cite{HS16}, we have shown that any pure jump subordinator of type-$(\alpha,\beta)$ can be decomposed into the summation of
two independent pure jump subordinators of type-$(\alpha,\beta)$ such that both of them satisfy (H) (see \cite[Theorem 4.2]{HS16}).

\begin{exa}\label{exa-4.14}
Let $0<\alpha_1<\beta_1<\alpha<\beta<1$. Suppose that $X_1$ is a pure jump subordinator of type-$(\alpha,\beta)$ satisfying
(H) and $X_2$ is a pure jump subordinator of type-$(\alpha_1,\beta_1)$ which is independent of $X_1$. We will prove
below that both $X_1+X_2$ and $X_1-X_2$ satisfy (H).

We denote by $\psi_1$ and $\psi_2$ the L\'{e}vy-Khintchine
exponents of $X_1$ and $X_2$, respectively.  Note that
$\overline{\psi_2}$ is the L\'{e}vy-Khintchine exponent $-X_2$. By
\cite[(4.5) and (4.6)]{HS16}, we find that there exist two
positive constants $c_1$ and $c_2$  such that
\begin{eqnarray}\label{exa-4.14-a}
1+{\rm Re}\psi_1(z)\geq 1+c_1|z|^{\alpha},\ \mbox{\rm for all}\
|z|\geq 1,
\end{eqnarray}
and
\begin{eqnarray*}
|\psi_2(z)|\leq c_2|z|^{\beta_1},\
\mbox{\rm for all}\ |z|\geq 1.
\end{eqnarray*}
Hence there exists a constant $c>0$ such that
$$
|\psi_2(z)|\leq c(1+{\rm Re}\psi_1(z)),\ \forall z\in \mathbf{R}.
$$
By (\ref{exa-4.14-a}) and \cite{HW42}, we know that $X_1$ has
transition densities and thus both $X_1+X_2$ and $X_1-X_2$ have transition densities.
Therefore, both $X_1+X_2$ and $X_1-X_2$ satisfy (H) by Theorem \ref{thm-4.1} and Proposition \ref{pro-3.12}.

\end{exa}

\bigskip

{ \noindent {\bf\large Acknowledgments} \vskip 0.1cm  \noindent We
acknowledge the support of NNSFC (Grant No. 11371191) and NSERC
(Grant No. 311945-2013).}


\begin{thebibliography}{1234}


\bibitem{B96} Bertoin J.: L\'{e}vy processes. Cambridge University Press, Cambridge (1996).

\bibitem{BG61} Blumenthal R.M., Getoor R.K.: Sample functions of stochastic processes with stationary independent increments. J. Math. Mech., \textbf{10}, 493-516 (1961).





\bibitem{Br71} Bretagnolle J.:  R\'{e}sults de Kesten sur les processus \`{a}
accroissements ind\'{e}pendants. S\'{e}minare de Probabilit\'{e}s
V, Lect. Notes in Math., Vol. 191, Springer-Verlag, Berlin, 21-36
(1971).

\bibitem{Fi14} Fitzsimmons P.J.: Gross's Bwownian motion fails to satisfy the polarity principle. Rev. Roumaine Math. Pures Appl. {\bf 59},  87-91 (2014).

%
%
%
\bibitem{F75}Forst
G.:  The definition of energy in non-symmetric translation
invariant Dirichlet spaces. Math. Ann. {\bf 216}, 165-172 (1975).

%



%
%



%
%
\bibitem{HN16} Hansen W., Netuka I.: Hunt's hypothesis (H) and triangle property of the Green function. Expo. Math. {\bf 34}, 95-100 (2016).

    \bibitem{HW42}  Hartman P., Wintner A.:   On the infinitesimal
generators of integral convolutions. Amer. J. Math. {\bf 64},
273-298 (1942).
%
%




\bibitem{HS12} Hu Z.-C., Sun W.: Hunt's hypothesis (H) and
Getoor's conjecture for L\'{e}vy processes. Stoch. Proc. Appl.
\textbf{122}, 2319-2328 (2012).

\bibitem{HS16} Hu Z.-C., Sun W.: Further study on Hunt's hypothesis (H) for L\'{e}vy processes.  Sci. China Math., {\bf 59}, 2205-2226 (2016).

\bibitem{HSZ15} Hu Z.-C., Sun W., Zhang J.: New results on Hunt's hypothesis (H)  for L\'{e}vy processes. Potential Anal., {\bf 42}, 585-605 (2015).




\bibitem{Ka76} Kanda M.:  Two theorems on capacity for Markov processes
with stationary independent increments. Z. Wahrsch. verw. Gebiete
{\bf 35}, 159-165 (1976).



\bibitem{Ke69} Kesten H.:  Hitting probabilities of single points for
processes with stationary  independent increments. Memoirs of the
American Mathematical Society, No. 93, American Mathematical
Society, Providence, R.I. (1969).


%
%






\bibitem{R88} Rao M.:  Hunt's hypothesis for L\'{e}vy
processes. Proc. Amer. Math. Soc. {\bf 104},  621-624 (1988).

\bibitem{Sa99} Sato K.: L\'{e}vy processes and infinitely divisible distributions. Cambridge University Press, Cambridge (1999).


\end{thebibliography}
\end{document}